\input amstex
\documentstyle{amsppt}
\magnification=1200

\def\C{\Bbb C}
\def\R{\Bbb R}
\def\Z{\Bbb Z}
\def\wtA{\widetilde A}
\def\wtB{\widetilde B}
\def\wtphi{\widetilde\varphi}
\def\wtpsi{\widetilde\psi}
\def\F{{_0}F_1}
\def\FF{{_1}F_1}
\def\al{\alpha}
\def\be{\beta}
\def\ga{\gamma}
\def\ka{\kappa}
\def\si{\sigma}
\def\Ga{\Gamma}
\def\X{\frak X}
\def\Y{\frak Y}
\def\bY{\overline{\Y}}
\def\D{\frak D}

\def\sh{\operatorname{sh}}
\def\ch{\operatorname{ch}}
\def\tht{\thetag}

\def\Kpp{\Cal K_{++}}
\def\Kmm{\Cal K_{--}}
\def\Kpm{\Cal K_{+-}}
\def\Kmp{\Cal K_{-+}}
\def\KK{\widehat{\Cal K}}
\def\kpp{k_{++}}
\def\kpm{k_{+-}}

\TagsOnRight
\NoRunningHeads

\topmatter
\title Point processes and \\
the infinite
symmetric group. \\
Part V: Analysis of the matrix Whittaker kernel
\endtitle
\author Grigori Olshanski \endauthor
\address{\rm  Dobrushin Mathematics Laboratory, 
Institute for Problems of Information Transmission, 
Bolshoy Karetny 19, 101447 Moscow GSP--4, Russia.} {\it
E-mail\/}: {\tt olsh\@iitp.ru, olsh\@glasnet.ru } 
\endaddress

\thanks Supported by the Russian Foundation for Basic Research 
under grant 98-01-00303 and by
the Russian Program for Support of Scientific Schools under grant
96-15-96060. 
\endthanks
\abstract
The matrix Whittaker kernel has been introduced by A.~Borodin in Part
IV of the present series of papers. This kernel describes a point
process --- a probability measure on a space of countable point
configurations. The kernel is expressed in terms of the Whittaker
confluent hypergeometric functions. It depends on two parameters and
determines a $J$-symmetric operator $K$ in 
$L^2(\Bbb R_+)\oplus L^2(\Bbb R_+)$. 

It turns out that the operator $K$ can be represented in the form
$L(1+L)^{-1}$, where $L$ is a rather simple integral operator: the
kernel of $L$ is expressed in terms of elementary functions only.
This is our main result; it elucidates the nature of the matrix
Whittaker kernel and makes it possible to directly verify the
existence of the associated point process. 

Next, we show that the matrix Whittaker kernel can be degenerated to
a family of kernels expressed through the Bessel and Macdonald
functions. In this way one can obtain both the well--known Bessel
kernel (which arises in random matrix theory) and certain
interesting new kernels. 
\endabstract
\toc
\widestnumber\head{5.}
\head {} Introduction \endhead
\head 1. A model \endhead 
\head 2. Main result \endhead
\head 3. Diagonalization of the kernels \endhead
\head 4. The matrix tail kernel \endhead
\head 5. Degeneration to a Bessel--type kernel \endhead
\head {} References \endhead
\endtoc

\endtopmatter
\document
\head Introduction \endhead

The present paper continues a series of papers by Alexei Borodin and
the author: \cite{P.I -- P.IV}. In \cite{P.VI} we give a summary
of the results. We refer to \cite{P.VI} for a detailed introduction to
the subject and motivation.

One of the main conclusions of our work can be stated as follows:
certain stochastic point processes, which originate in harmonic
analysis on the infinite symmetric group, turn out to be close to
point processes arising in scaling limit of certain random matrix
ensembles. 

The common feature of the point processes in question is that their
correlation functions are given by determinantal expressions
$$
\rho_n(x_1,\dots,x_n)=\det[K(x_i,x_j)],  \tag0.1
$$
where $K(x,y)$ is a kernel in two real variables. In principle, all
the characteristics of such processes can be extracted from the
corresponding kernels, though in practice this often requires
a lot of work. 

Random matrix theory produces a variety of interesting kernels. Among
them are the Bessel kernel and the Airy kernel, which are expressed
through the Bessel functions $J_\nu$ and the Airy function,
respectively. About these kernels, see \cite{F, NS, NW, TW1--3}. 

Our work leads to a new family of kernels, which are expressed
through the Whittaker functions \footnote{The Whittaker functions are
certain confluent hypergeometric functions. They are eigenfunctions of
a second order differential operator on the semiaxis $x>0$ and have
exponential decay as $x\to+\infty$.} -- the {\it Whitaker
kernel\/} and the {\it matrix Whittaker kernel,\/} both depending on
two real parameters.

The Whitaker kernel has been introduced in \cite{P.II} and studied in
detail in \cite{P.III}. It describes a point process on 
$\R_+$ whose `particles' are accumulated about zero. 

The matrix Whittaker kernel, which is the object of study in the
present paper, has been introduced in \cite{P.IV}. It describes a larger
process, which lives on $\R\setminus\{0\}=\R_+\cup\R_-$. This kernel
is conveniently written as a $2\times 2$ matrix whose entries are
kernels (or integral operators) on $\R_+$, which explains the term
`matrix kernel'. 

The matrix Whittaker kernel defines an operator in 
$L^2(\R_+)\oplus L^2(\R_+)$ which turns out to be $J$-symmetric,
i.e., symmetric with respect to an indefinite inner product
$[f,g]=(Jf,g)$; specifically, 
$J=\bmatrix \operatorname{id} & 0\\ 0 & -\operatorname{id}
\endbmatrix$. Perhaps, this is the first example of a
$J$-symmetric kernel occuring in formula \tht{0.1} for the correlation
functions.  

The purpose of the present paper is to try to elucidate the nature of
the matrix Whittaker kernel and its relationship to other kernels.

In section 1, we explain some general properties of the point
processes governed by $J$-symmetric kernels on a simple model (finite
state space). This enables us to avoid unnecessary technicalities in
questions which are essentially based on elementary linear algebra.

In section 2, we employ the formulas of section 1 as a prompt to
derive the following result about the structure of the matrix 
Whittaker kernel (Theorem 2.4). Let $K$ denote the $J$-symmetric
operator in $L^2(\R_+)\oplus L^2(\R_+)$ mentioned above. We show that
the operator $L=K(1-K)^{-1}$ has the form
$$
L=\bmatrix 0 & A \\ -A' & 0 \endbmatrix\,, \tag0.2
$$
where $A$ is a real integral operator on $\R_+$ and $A'$ is the
transposed operator. It is worth noting that the kernel of $A$ is
given by a simple expression involving no special functions at all;
the Whittaker functions arise when we pass from $L$ to $K$. 

The passage from $K$ to $L$ is not a pure formal trick. As is
explained in section 1, in a simpler situation of processes with
finite point configurations (which is not the case of our process on
$\R_+\cup\R_-$), the operator $L$ has a clear meaning: its kernel
describes the distribution functions of the process in the whole
state space. In such a situation, vanishing of the diagonal blocks of
$L$ also admits a nice interpretation (see Proposition 1.7).  

For our process, the point configurations are almost surely infinite,
and the basic formula expressing the distribution functions in the whole
state space through the operator $L$ (see \tht{1.1}) becomes
incorrect. However, it is tempting to combine it with the
simple expression for the operator $L$ given in Theorem 2.4 to derive
meaningful conclusions about our process. In this direction, we
(Borodin and I) have some conjectures but no rigorous results. 

Section 3 is devoted to the spectral analysis of the matrix Whittaker
kernel $K(x,y)$. We exhibit a continual basis in 
$L^2(\R_+)\oplus L^2(\R_+)$ diagonalizing the corresponding operator
$K$. The elements of this basis are (within a simple factor) certain
Whittaker functions; they are eigenfunctions of a Sturm--Liouville
differential operator. This result provides additional information
about the nature of the matrix Whittaker kernel $K(x,y)$. It also
provides a way to verify directly that $K(x,y)$ obeys certain
conditions that ensure the existence of a point process with the
correlation functions \tht{0.1}. 

In section 4 we deal with the `tail process'. This
point process is a stationary process on $\R\cup\R$. It
describes (via an appropriate scaling limit) the asymptotical behavior of
our initial random configurations on $\R_+\cup\R_-$ near zero. This
extends the results of sections 3--4 of \cite{P.III}.

In section 5 we study the matrix Whittaker kernel $K(x,y)$ from
another point of view: hierarchy of special functions of
hypergeometric type. It is well--known that the Whittaker functions
can be degenerated to the Bessel functions. Employing this fact, we
compute a scaling limit of the kernel $K(x,y)$. The limit transition
is rather curious: we let one of the two parameters of the kernel,
$a$, tend to infinity inside the set of the form $a_0+2\Z$ with
$a_0\in\R$ fixed. As a result, we get a (still two--parametric) family
of matrix kernels expressed through the Bessel functions of the first
kind $J_\nu$ and the modified Bessel functions of the third kind
(Macdonald functions) $K_\nu$. Taking the diagonal blocks, we get
two sorts of kernels on $\R_+$: one (expressed through the functions
$J_\nu$) is a slight generalization of the conventional Bessel kernel
\cite{F, NS, NW, TW2}, and another (expressed through the Macdonald functions
$K_\nu$) is, perhaps, a new example; we called it the Macdonald
kernel. 

It would be interesting to understand whether the Macdonald kernel is
somehow related to random matrix ensembles.

I am grateful to Alexei Borodin for numerous discussions and to Craig A.
Tracy for drawing my attention to the papers \cite{MTW}, \cite{T}, \cite{TW4}.
 
\head 1. A model  \endhead

In this section we fix a finite set $\X$ which will serve as a ``state
space''. We shall deal with kernels $K(x,y)$, $L(x,y)$ on $\X\times\X$
which will also be considered as matrices of order $|\X|$ with the
rows and columns indexed by the elements of $\X$. The
kernels can be real or complex. We take the counting measure on $\X$ and
form the corresponding (finite--dimensional) Hilbert space $L^2(\X)$.
Any kernel defines an operator in $L^2(\X)$; we shall denote it by
the same letter as the kernel. 

Let $\Xi$ stand for the space of configurations in $\X$ with no
multiple points. Thus, a configuration is simply a (possibly, empty)
subset $\xi\subseteq\X$, and the total number of configurations is
equal to $2^{|\X|}$. 

Let $L(x,y)$ be a kernel on $\X$. For a configuration $\xi$ we shall
denote by $L_\xi$ the submatrix in $L$ formed by the rows and columns
from $\xi$, so that $\det L_\xi$ is a principal minor of $L$. We
agree that $\det L_{\varnothing}=1$.

\proclaim{Proposition 1.1} Assume $L$ is a real or complex kernel on
$\X$ such that all its principal minors $\det L_\xi$ are real and
nonnegative. Then there exists a probability distribution on $\Xi$
with the weights
$$
\operatorname{Prob}\{\xi\}=\frac{\det L_\xi}{\det(1+L)}\,. \tag1.1
$$
\endproclaim

Note that the assumptions on $L$ imply that the matrix $1+L$ is
invertible.

\demo{Proof} This is easy. \qed
\enddemo

We shall consider this probability distribution as a point process
on the (finite) state space $\X$. The next step is to look at the
correlation functions.

\proclaim{Proposition 1.2} The correlation functions of the above
process are given by the determinantal formula
$$
\rho_n(x_1,\dots,x_n)=\det[K(x_i,x_j)]_{1\le i,j\le n}\,, \tag1.2
$$
where $n=1,2,\dots|\X|$, the points $x_1,\dots,x_n$ are pairwise
distinct, and $K$ is given by $K=L(1+L)^{-1}$.
\endproclaim

\demo{Proof} An elegant proof based on the generating functional
of the process is given in \cite{DVJ}, Exercise 5.4.7. \qed
\enddemo

We shall always assume that $K$ and $L$ are related to each other by
the transformations
$$
K=L(1+L)^{-1}, \quad L=K(1-K)^{-1},\tag1.3
$$
which are well--defined provided that the matrices $1+L$ and $1-K$
are invertible. 

\proclaim{Proposition 1.3} Assume $L$ is Hermitian nonnegative:
$L=L^*\ge0$. Then $L$ satisfies the assumptions of Proposition 1.1 and
$K$ is Hermitian satisfying $0\le K<1$. Conversely, if $K$ is
Hermitian and $0\le K<1$ then the corresponding $L$ exists and is
Hermitian nonnegative.
\endproclaim

\demo{Proof} This follows from \tht{1.3}. \qed
\enddemo

Thus, we dispose of a family of (finite) point processes governed by
the Hermitian, nonnegative, strictly contractive kernels $K$.
Or, equivalently, by the Hermitian nonnegative kernels $L$.

Henceforth we fix a partition of $\X$ into disjoint union of two
subsets:
$$
\X=\X_1\cup\X_2,
$$ 
we decompose the Hilbert space $L^2(\X)$ into the direct sum
$$
L^2(\X)=L^2(\X_1)\oplus L^2(\X_2),\tag1.4
$$
and we write any kernel (matrix) $A$ on $\X$ in the block form
$$
A=\bmatrix A_{11} & A_{12}\\ A_{21} & A_{22} \endbmatrix\,, 
$$
where $A_{ij}$ acts from $L^2(\X_j)$ to $L^2(\X_i)$.

We endow the space $L^2(\X)$ with the indefinite inner product
determined by the matrix $\bmatrix 1 & 0\\ 0 & -1 \endbmatrix$.  A
kernel $A$ is called {\it $J$-Hermitian\/} if the corresponding
operator is Hermitian with respect to the indefinite inner product.
Equivalently, in terms of blocks, 
$$
(A_{11})^*=A_{11},\quad (A_{22})^*=A_{22},\quad (A_{12})^*=-A_{21}. \tag1.5
$$

Note that if one of the kernels $K$, $L$ is $J$-Hermitian, and the
correspondence \tht{1.3} makes sense, then another kernel is
$J$-Hermitian, too.

\proclaim{Proposition 1.4} Assume that $L$ is $J$-Hermitian and the
diagonal blocks $L_{11}$, $L_{22}$ are nonnegative. Then all the
principal minors $\det L_\xi$ are real nonnegative.
\endproclaim

\demo{Proof} Replacing $L$ by $L+\varepsilon1$, where
$\varepsilon>0$ is arbitrary, we may assume that the diagonal blocks
$L_{11}$, $L_{22}$ are strictly positive. We shall prove that under this
assumption, all the principal minors $\det L_\xi$ are strictly
positive. Next, replacing $L$ by $L_\xi$, we remark that it suffices
to prove this claim for $\det L$ only. 

Thus, we have to prove that $\det L>0$ provided that $L$ is
$J$-Hermitian and $L_{11}>0$, $L_{22}>0$. We employ the well--known
formula for the determinant of a block matrix,
$$
\det L=\det L_{11}\cdot 
\det(L_{22}-L_{21}L_{11}^{-1}L_{12})\,
$$
which makes sense because $L_{11}$ is invertible. We have
$$
L_{11}>0, \quad 
L_{22}-L_{21}L_{11}^{-1}L_{12}=L_{22}+L_{12}^*L_{11}^{-1}L_{12}>0\,,
$$
which implies that the both determinants are strictly positive. \qed
\enddemo

\proclaim{Proposition 1.5} In terms of blocks, the correspondence
$L\leftrightarrow K$ takes the following form.

The transform $L\mapsto K$:
$$
\gather
K_{11}=(L_{11}-L_{12}(1+L_{22})^{-1}L_{21})
(1+L_{11}-L_{12}(1+L_{22})^{-1}L_{21})^{-1}\tag1.6a\\
K_{22}=(L_{22}-L_{21}(1+L_{11})^{-1}L_{12})
(1+L_{22}-L_{21}(1+L_{11})^{-1}L_{12})^{-1}\tag1.6b\\
K_{12}=(1-K_{11})L_{12}(1+L_{22})^{-1}\tag1.6c\\
K_{21}=(1-K_{22})L_{21}(1+L_{11})^{-1}\,. \tag1.6d
\endgather
$$

The inverse transform $K\mapsto L$:
$$
\gather
L_{11}=(K_{11}-K_{12}(1+K_{22})^{-1}K_{21})
(1+K_{11}-K_{12}(1+K_{22})^{-1}K_{21})^{-1}\tag1.7a\\
L_{22}=(K_{22}-K_{21}(1+K_{11})^{-1}K_{12})
(1+K_{22}-K_{21}(1+K_{11})^{-1}K_{12})^{-1}\tag1.7b\\
L_{12}=(1+L_{11})K_{12}(1-K_{22})^{-1}\tag1.7c\\
L_{21}=(1+L_{22})K_{21}(1-K_{11})^{-1}\,.\tag1.7d
\endgather
$$

Here we assume that all the necessary inverse matrices exist.
\endproclaim

\demo{Proof} Write  the
equality $K=L(1+L)^{-1}$ or, equivalently, $K(1+L)=L$ in terms of
blocks:
$$
\gather
K_{11}(1+L_{11})+K_{12}L_{21}=L_{11}\\
K_{11}L_{12}+K_{12}(1+L_{22})=L_{12}\\
K_{22}(1+L_{22})+K_{21}L_{12}=L_{22}\\
K_{22}L_{21}+K_{21}(1+L_{11})=L_{21}\,.
\endgather
$$
{}From this system one readily gets the relations \tht{1.6a-d}.

The second group of relations, \tht{1.7a-d}, is verified in exactly
the same way. \qed
\enddemo

\proclaim{Proposition 1.6} The transforms $L\mapsto K$ and 
$K\mapsto L$ define a bijective correspondence between

{\rm(i)} the $J$-Hermitian kernels $L$ such that $L_{11}\ge0$,
$L_{22} \ge0$, \par
\noindent and  

{\rm(ii)} the $J$-Hermitian kernels $K$ such that $K_{11}<1$,
$K_{22}<1$, 
$$
K_{11}+K_{12}(1-K_{22})^{-1}K_{21}\ge0, \quad
K_{22}+K_{21}(1-K_{11})^{-1}K_{12}\ge0 \tag1.8
$$
(note that the latter two inequalities are stronger than $K_{11}\ge0$,
$K_{22}\ge0$). 

For these two sets of kernels all the inverse matrices in the
formulas of Proposition 1.5 exist. 
\endproclaim

\demo{Proof} We have $1+L_{22}>0$, so that $1+L_{22}$ is invertible.
Further, the matrix
$$
P_{11}:=L_{11}-L_{12}(1+L_{22})^{-1}L_{21}=
L_{11}+L_{12}(1+L_{22})^{-1}L_{12}^*
$$
is nonnegative. Since the expression \tht{1.6a} for $K_{11}$ is
equivalent to 
$$
K_{11}=P_{11}(1+P_{11})^{-1},
$$
it makes sense and, moreover, $0\le K_{11}<1$. Likewise, the
expression \tht{1.6b} for $K_{22}$ also makes sense and we have 
$0\le K_{22}<1$. Consequently, the expressions \tht{1.6c-d} for $K_{12}$ and
$K_{21}$ also make sense. Thus, the transform $L\mapsto K$ is
well--defined. 

Now look at the relation \tht{1.7a}. It can be written as
$$
L_{11}=Q_{11}(1-Q_{11})^{-1},
$$
where 
$$
Q_{11}=K_{11}+K_{12}(1-K_{22})^{-1}K_{21}.
$$
We have
$$
1-Q_{11}=1-K_{11}-K_{12}(1-K_{22})^{-1}K_{21}=
(1-K_{11})+K_{12}(1-K_{22})^{-1}K_{12}^*>0,
$$
because $K_{11}<1$. Since $L_{11}\ge0$, we conclude $Q_{11}\ge0$,
which is the first inequality in \tht{1.8}. It is stronger than
$K_{11}\ge0$, because $Q_{11}\ge0$ means
$$
K_{11}-(K_{12}(1-K_{22})^{-1}K_{12}^*)\ge0.
$$
Likewise, we establish the second inequality in \tht{1.8}.

Thus, the transform $L\mapsto K$ sends any matrix $L$ satisfying (i)
into a matrix $K$ satisfying (ii). The inverse transform is justified
in the same way. \qed
\enddemo

Proposition 1.6 is a generalization of the trivial Proposition 1.3.
Together with Proposition 1.4 it yields a class of finite point
processes for which both the kernels $L$ and $K$ are given
explicitly. We recall that $L$ describes the distribution functions 
(see \tht{1.1}), while $K$ describes the correlation functions (see
\tht{1.2}).  

Now we shall impose more special conditions on $L$: we shall assume
that $L_{11}=0$, $L_{22}=0$. The meaning of this restriction is
demonstrated by the following result.

\proclaim{Proposition 1.7} Consider the process governed by a
$J$-Hermitian kernel $L$ with $L_{11}\ge0$, $L_{22}\ge0$. Then the
condition $L_{11}=0$, $L_{22}=0$ exactly means that the process is
concentrated on the configurations $\xi\subset\X$ with the property 
$$
|\xi\cap\X_1|=|\xi\cap\X_2|.\tag1.9
$$
\endproclaim

\demo{Proof} Assume $L_{11}=0$, $L_{22}=0$, and let $\xi$ be a nonempty
configuration. Its weight \tht{1.1} is proportional to the value of the
principal minor $\det L_\xi$. Remark that $L_\xi$ is a block matrix
whose diagonal blocks are zero. Such a matrix can be nondegenerate
only if the blocks are of the same size, which means \tht{1.9}.

Conversely, let \tht{1.9} hold. The diagonal blocks $L_{11}$, $L_{22}$
are nonnegative, so that if one of them is nonzero, then it has a
nonzero diagonal entry. That is, there exists a point $x\in\X$ such
that $L(x,x)>0$. But then the one--point configuration $\xi=\{x\}$
has a nonzero weight which contradicts to the assumption \tht{1.9}. This
contradiction implies that the diagonal blocks must be zero. \qed
\enddemo

\proclaim{Proposition 1.8} The transforms $L\mapsto K$ and $K\mapsto L$
define a bijective correspondence between

{\rm(i)} the kernels $L$ of the form
$$
L=\bmatrix 0 & A\\ -B & 0\endbmatrix\,, \tag1.10
$$
where the matrix $1+AB$ is invertible (equivalently, $1+BA$ is
invertible) 

\noindent and

{\rm(ii)} the kernels $K$ of the form
$$
K=\bmatrix CD & C\\ DCD-D & DC\endbmatrix \,, \tag1.11
$$
where $1-CD$ is invertible (equivalently, $1-DC$ is invertible).

In terms of the blocks, this correspondence takes the form
$$
\gather
C=(1+AB)^{-1}A=A(1+BA)^{-1}, \quad D=B,\tag1.12a\\
A=C(1-DC)^{-1})=(1-CD)^{-1}, \quad B=D.\tag1.12b
\endgather
$$

In particular, 
$$
1-CD=(1+AB)^{-1}, \quad 1-DC=(1+BA)^{-1}.\tag1.13
$$
\endproclaim

\demo{Proof} This is a direct consequence of the formulas of
Proposition 1.5 and the identity
$$
X(1\pm YX)^{-1}=(1\pm XY)^{-1}X. \tag1.14
$$

Indeed, take the kernel $L$ of the form \tht{1.10}. Applying
\tht{1.6a-d} we get 
$$
K=\bmatrix AB(1+AB)^{-1} & (1+AB)^{-1}A \\
-(1+BA)^{-1}B & BA(1+BA)^{-1} \endbmatrix\,. \tag1.15
$$
Using the identity \tht{1.14 } we verify that $K$ has the desired form
with $C,D$ as indicated in \tht{1.12a}.

Conversely, starting with the kernel $K$ of the form \tht{1.11} and
applying \tht{1.7a-d} together with the identity \tht{1.14} we
readily verify that $L$ has the form \tht{1.10} with $A,B$ indicated
in \tht{1.12b}. \qed 
\enddemo

The following corollary will serve as a prompt for the main result of
the next section (Theorem 2.4).

\proclaim{Corollary 1.9} Let $K$ be a $J$-Hermitian kernel of the
form \tht{1.11}. Then 
$$
L=\bmatrix 0 & D^* \\ -D & 0 \endbmatrix\,.
$$
\endproclaim

\demo{Proof} By Proposition 1.8, $L$ is given by the formula
\tht{1.10} with $B=D$. Since $K$ is $J$-Hermitian, $L$ is
$J$-Hermitian, too. This implies $A=B^*=D^*$. \qed
\enddemo

\example{Remark 1.10} Let us return to the very beginning of the
section, where we introduced a class of point processes in $\X$
governed by kernels $L$ with 
nonnegative principal minors (formula \tht{1.1}). Let $\Y\subseteq\X$ be
a subset. Given a point process in $\X$ we can define its
`truncation', which is a point process in $\Y$: the latter process is
the image of the former under the map $\xi\mapsto\xi\cap\Y$. Clearly,
the correlation functions of the truncated process are obtained
simply by restricting the correlation functions of the initial
process. So, if the initial process belongs to our class (i.e.,
is given by the formula \tht{1.1}) then the truncated process also
belongs to this class, and the `truncated' kernel $K$ is obtained
simply by restricting the initial kernel $K$ to $\Y\times\Y$. The
corresponding transformation of the kernel $L$ is more complicated.
To describe it, write the kernel $L$ in the block form with respect
to the partition $\X=\Y\cup\bY$:
$$
L=\bmatrix L_{\Y\Y} & L_{\Y\bY}\\ L_{\bY\Y} & L_{\bY\bY}\endbmatrix
$$
Then the transformed kernel is equal to
$$
L_{\Y\Y}-L_{\Y\bY}(1+L_{\bY\bY})^{-1}L_{\bY\Y}.
$$
This formula follows from \tht{1.6a}.
\endexample

\head 2. Main result \endhead

Here we shall try to apply the general results of section 1 to the
matrix Whittaker kernel introduced in \cite{P.IV}. As we shall deal
with continual kernels instead of finite matrices, we shall need to
justify certain steps.  

We take as $\X$ the punctured real line $\R^*=\R\setminus\{0\}$ equipped
with Lebesgue measure $dx$. We write kernels on $\R^*$ in the block
form with respect to the partition 
$$
\R^*=\R_+\cup\R_-\,.
$$

We identify $\R_-$ with $\R_+$ via the map $x\mapsto |x|$, $x<0$,
which enables us to interpret each block as a kernel on $\R_+$ and
write
$$
L^2(\R^*)=L^2(\R_+)\oplus L^2(\R_+). \tag2.1
$$

Given a kernel, we shall denote the corresponding integral operator
by the same letter and we shall regard it as a Hilbert space operator.

We shall introduce the basic notation and then we shall recall
the definition of the matrix Whittaker kernel from \cite{P.IV}.

This kernel depends on two parameters $z,z'$ which play  symmetric
roles. The parameters satisfy the following restrictions:
$$
 \gathered 
\text{either $z,z'\in\C\setminus\Z$ and $z'=\bar z$,}\\
\text{or $z,z'\in\R\setminus\Z$ and $m<z,z'<m+1$ for a certain
$m\in\Z$.}
\endgathered\tag2.2
$$

Instead of $z,z'$ one can take the parameters $a,\mu$ defined by
$$
a=\frac{z+z'}2,\quad \mu=\frac{z-z'}2\,.  \tag2.3
$$
Then the conditions \tht{2.2} take the following form:
$$
\gathered 
\text{$a$ is always real and is not an integer when $\mu=0$;}\\
\text{$\mu$ is either pure imaginary or real; in the latter case}\\ 
\text{there exists $m\in\Z$ such that
$m+|\mu|<a<m+1-|\mu|$.}
\endgathered
\tag2.4
$$
This implies that if $\mu$ is real then $|\mu|<\frac12$.

The Whittaker function $W_{\ka,\mu}(x)$, $x>0$, with indices
$\ka,\mu$, can be defined in terms of the confluent hypergeometric
function $\FF$ as follows: 
$$
\aligned
x^{-1/2}e^{x/2}W_{\ka,\mu}(x) &=
\frac{\Ga(-2\mu)x^\mu}{\Ga(\frac12-\ka-\mu)}\,
\FF(\tfrac12-\ka+\mu;2\mu+1;x)\\
&+\frac{\Ga(2\mu)x^{-\mu}}{\Ga(\frac12-\ka+\mu)}\,
\FF(\tfrac12-\ka-\mu;-2\mu+1;x).
\endaligned \tag2.5
$$
This expression makes sense when $\mu\ne0$; when $\mu=0$ (the
so--called logarithmic case), it can be defined by a limit
transition. We have
$$
W_{\ka,\mu}=W_{\ka,-\mu}\,. \tag2.6
$$
We shall always deal with a real $\ka$ and a real or pure imaginary
$\mu$; then the Whittaker function takes real values.

The Whittaker function can be characterized  as
the only solution of the second order differential equation
$$
\frac{d^2W}{dx^2}+\left(-\frac14+\frac{\ka}x+
\frac{\frac14-\mu^2}{x^2}\right)W=0, \quad x>0, \tag2.7
$$
which has exponential decay at $+\infty$:
$$
W_{\ka,\mu}(x)=x^\ka e^{-x/2}\left(1+O\left(\frac1x\right)\right),
\quad x\to +\infty.  \tag2.8
$$
About the Whittaker function, see, e.g., \cite{E1}.

Let us fix $z,z'$ (equivalently, $a,\mu$) and abbreviate
$$
\gather
\varphi(x)=x^{-1/2}e^{x/2}W_{a+\frac12,\mu}(x),\tag2.9a \\
\varphi_-(x)=x^{-1/2}e^{x/2}W_{a-\frac12,\mu}(x),\tag2.9b \\
\psi(x)=x^{-1/2}e^{x/2}W_{-a+\frac12,\mu}(x),\tag2.9c \\
\psi_-(x)=x^{-1/2}e^{x/2}W_{-a-\frac12,\mu}(x).\tag2.9d \\
\endgather
$$

Finally, let 
$$
\si=\sqrt{\sin(\pi z)\sin(\pi z')}=
\sqrt{\frac{\cos(2\pi\mu)-\cos(2\pi a)}2} \tag2.10
$$ 
and note that, under our restrictions on the parameters, $\si^2$ is
always real and strictly positive; so, $\si$ is always real and
nonzero.

The {\it matrix Whittaker kernel\/}, as defined in \cite{P.IV}, has the
form 
$$
[K]=\bmatrix K_{++} & K_{+-} \\ K_{-+} & K_{--} \endbmatrix\,, \tag2.11
$$ 
where the blocks are the following kernels on $\R_+$:
$$
\gather
K_{++}(x,y)=
\frac1{\Ga(z)\Ga(z')}\,
\frac{\varphi(x)\varphi_-(y)-\varphi_-(x)\varphi(y)}{x-y} \tag2.12a\\
K_{--}(x,y)=
\frac1{\Ga(-z)\Ga(-z')}\,
\frac{\psi(x)\psi_-(y)-\psi_-(x)\psi(y)}{x-y}  \tag2.12b\\
K_{+-}(x,y)=
\frac{\si}{\pi}\,
\frac{\varphi(x)\psi(y)+zz'\varphi_-(x)\psi_-(y)}{x+y}  \tag2.12c\\
K_{-+}(x,y)=-K_{+-}(y,x)  \tag2.12d
\endgather
$$

Note that all the blocks are real kernels and the diagonal kernels
$K_{++}$, $K_{--}$ are symmetric; together with the last relation
this implies that $[K]$ is $J$-symmetric.

\proclaim{Proposition 2.1} The matrix Whittaker kernel $[K]$ can be
written in the form
$$
[K]=\bmatrix CD & C\\ DCD-D & DC \endbmatrix,  \tag2.13
$$
where 
$$
\gather
C(x,y)=K_{+-}(x,y),   \tag2.14\\
D(x,y)=\frac\si\pi\,\left(\frac xy\right)^{-a}\,
\frac{e^{-\frac{x+y}2}}{x+y}\,,   \tag2.15
\endgather
$$
and, by definition, the products $CD$, $DC$, $DCD$ are superpositions of
integral operators: 
$$
\gathered
(CD)(x,y)=\int_0^{+\infty}C(x,s)D(s,y)ds\\
(DC)(x,y)=\int_0^{+\infty}D(x,s)C(s,y)ds\\
(DCD)(x,y)=\int_0^{+\infty}\int_0^{+\infty}
D(x,s_1)C(s_1,s_2)D(s_2,y)ds_1ds_2\,.
\endgathered
\tag2.16
$$
\endproclaim

\demo{Proof} This is merely a reformulation of the results of
\cite{P.IV}, section 2. Indeed, let 
$$
N(x,y)=\frac\si\pi\, \left(\frac xy\right)^a \,e^{-\frac{x+y}2}\,
K_{+-}(x,y)\,,
\quad W(x,y)=\frac1{x+y}\,,  \tag2.17
$$
and let $R_1$ and $R_2$ be the following
multiplication operators:
$$
(R_1f)(x)=\pi\si^{-1}x^{-a}e^{x/2}f(x),\quad
(R_2)f(x)=x^{-a}e^{-x/2}f(x). \tag2.18
$$
In \cite{P.IV}, section 2, it was proved that
$$
K=\bmatrix R_1 & 0 \\ 0 & R_2 \endbmatrix \cdot 
\bmatrix NW & N \\ WNW-W & WN \endbmatrix\cdot 
\bmatrix R_1^{-1} & 0 \\ 0 & R_2^{-1} \endbmatrix\,,
\tag2.19
$$
which means, in particular, that all the integrals involved in
\tht{2.16} make sense. It follows that $[K]$ has the form \tht{2.13}
with 
$$
C=R_1NR_2^{-1},\quad D=R_2WR_1^{-1},
$$ 
which are exactly \tht{2.14}, \tht{2.15}. \qed
\enddemo

\proclaim{Proposition 2.2} Assume $-1/2<a<1/2$. Then the integral
operator $D$ with the kernel $D(x,y)$ as defined in \tht{2.15} is
bounded in $L^2(\R_+,dx)$. 
\endproclaim

\demo{Proof} Set
$$
D'(x,y)=\left(\frac xy\right)^{-a}\frac1{x+y}\,.  \tag2.20
$$
We have $D=(\si/\pi)MD'M$, where $M$ stands for the operator of
multiplication by the bounded function $e^{-x/2}$, so that it
suffices to check that $D'$ is bounded.

Passing to new variables $\xi,\eta$ such that $x=e^{-\xi}$,
$y=e^{-\eta}$, we transform $D'$ to the integral operator in
$L^2(\R,d\xi)$ with the kernel
$$
\Cal D'(\xi,\eta)=\frac{e^{-a(\xi-\eta)}}
{e^{(\xi-\eta)/2}+e^{-(\xi-\eta)/2}}\,.  \tag2.21
$$
Since this kernel is translation invariant, the Fourier transform
takes $\Cal D'$ to a multiplication operator. It remains to check
that the latter is bounded.

To this end we employ the formula
$$
\gathered
Fourier_u\left\{\frac{e^{\al v}}{e^{\be v}+e^{-\be v}}\right\}
=\frac\pi{2\be}\,\frac1{\cos(\pi(-\frac\al{2\be}+\frac{iu}{2\be}))}\\
\be>0,\quad \Re\al<\be
\endgathered
\tag2.22
$$
where
$$
Fourier_u\{f(v)\}=g(u)\quad\text{means}\quad
g(u)=\int e^{iuv}f(v)dv,
$$
see \cite{E2, 3.2(15)}.

Applying \tht{2.22} we get that after the Fourier transform the
operator $\Cal D'$ becomes the operator of multiplication by the
function 
$$
\frac{\pi}{\cos(\pi(a+iu))}\,,
$$
which is bounded by the assumption $|a|<1/2$. \qed
\enddemo

\proclaim{Proposition 2.3} Assume $-1/2<a<1/2$. Then the integral
operator $D$ with the kernel $C(x,y)=K_{+-}(x,y)$ as defined by the
formulas \tht{2.14}, \tht{2.12c}, \tht{2.10}, \tht{2.9}, \tht{2.5},
is bounded in $L^2(\R_+,dx)$.  
\endproclaim

\demo{Proof} Assume first that $\mu\ne0$. Then, according to \tht{2.5},
we can write
$$
x^{-1/2}W_{\ka,\mu}(x)=x^\mu(\dots)+x^{-\mu}(\dots),  \tag2.23
$$
where each of the two expressions denoted by dots is equal to a confluent
hypergeometric function multiplied by an exponential factor and so is
analytic on the whole real axis. Consequently, the kernel $C(x,y)$
can be represented in the form 
$$
C(x,y)=\frac{x^\mu y^\mu a_{++}(x,y)
+x^{-\mu}y^{-\mu} a_{--}(x,y)
+x^\mu y^{-\mu} a_{+-}(x,y)
+x^{-\mu}y^\mu a_{-+}(x,y)}
{x+y}\,,
$$
where the four functions $a_{++}$, $a_{--}$, $a_{+-}$, $a_{-+}$ are
analytic about $(0,0)$. 

The key observation is that
$$
a_{++}(0,0)=a_{--}(0,0)=0.
$$
Indeed, it follows from \tht{2.5} that
$$
\gather
a_{++}(0,0)=\frac\si\pi\,(\Ga(-2\mu))^2\,
\left(\frac1{\Ga(-z)\Ga(z')}+\frac{zz'}{\Ga(-z+1)\Ga(z'+1)}\right)\\
a_{--}(0,0)=\frac\si\pi\,(\Ga(-2\mu))^2\,
\left(\frac1{\Ga(-z')\Ga(z)}+\frac{zz'}{\Ga(-z'+1)\Ga(z+1)}\right)\,,
\endgather
$$
and the both expressions in the parentheses are equal to zero.

Consider the kernel
$$
C_0(x,y)=\frac{x^\mu y^{-\mu}a_{+-}(0,0)+x^{-\mu}y^\mu a_{-+}(0,0)}
{x+y}\,,
$$
and set
$$
C'(x,y)=C(x,y)-\chi_{[0,1]}(x)C_0(x,y)\chi_{[0,1]}(y),
$$
where $\chi_{[0,1]}$ stands for the characteristic function of
$[0,1]$. The boundedness of the integral operator $C$ will follow
from the two claims: 
\medskip

$\bullet$ {\it Claim\/} 1. The kernel $C_0(x,y)$ defines a bounded
operator in $L^2(\R_+,dx)$.

$\bullet$ {\it Claim\/} 2. The function $C'(x,y)$ is square integrable
on $\R_+^2$.
\medskip

Indeed, $C(x,y)$ is the sum of $\chi_{[0,1]}(x)C_0(x,y)\chi_{[0,1]}(y)$
and $C'(x,y)$. Since the function $\chi_{[0,1]}$ is bounded, the
first summand defines a bounded operator by Claim 
1, and the second summand defines a Hilbert--Schmidt (hence bounded)
operator by Claim 2.  

Let us check Claim 1. It suffices to show that the kernel 
$$
\left(\frac xy\right)^{\pm\mu}\,\frac1{x+y}
$$
defines a bounded operator. Arguing as in the proof of Proposition
2.2 we reduce this to the boundedness of the function
$$
\frac\pi{\cos(\mp\pi\mu+i\pi u)}\,.
$$
Since $\mu\in i\R$ or $-1/2<\mu<1/2$ this is obvious.

Let us check Claim 2. It suffices to show that $C'$ is square
integrable both in the square $[0,1]^2$ and in its complement
$\R_+^2\setminus [0,1]^2$. 

In the second region, $C'$ coincides with
$C$. Recall that $C(x,y)=K_{+-}(x,y)$ (see \tht{2.14}) and look at
the expression \tht{2.12c} for $K_{+-}(x,y)$. Outside 
$[0,1]^2$, the fraction $\frac1{x+y}$ is bounded, so that it remains to
prove the square integrability of the numerator of \tht{2.12c}. To 
this end it suffices to prove that each of the functions
$\varphi(x)\psi(y)$, $\varphi_-(x)\psi_-(y)$ is square integrable in
our region. Actually we can claim that they are square integrable in
the whole quadrant $\R_+^2$. Indeed, this reduces to the fact that
each of the four functions $\varphi$, $\psi$, $\varphi_-$, $\psi_-$
is square integrable in $\R_+$. Each of these functions is of the form
$x^{-1/2}W_{\ka,\mu}(x)$. The latter function has exponential decay
at infinity, hence it is square integrable near infinity. Finally,
near zero, it behaves as indicated in \tht{2.23}. Since $\mu$ is either
pure imaginary or satisfies $-1/2<\mu<1/2$ we get square
integrability about zero as well. 

Let us examine the kernel $C'(x,y)$ in the square $[0,1]^2$; here it
coincides with $C(x,y)-C_0(x,y)$. By the definition of $C_0$ we have
$$
C'(x,y)=\frac{x^\mu y^\mu b_{++}(x,y)
+x^{-\mu}y^{-\mu} b_{--}(x,y)
+x^\mu y^{-\mu} b_{+-}(x,y)
+x^{-\mu}y^\mu b_{-+}(x,y)}
{x+y}\,,
$$
where the functions 
$$
\gather
b_{++}(x,y)=a_{++}(x,y),\quad 
b_{--}(x,y)=a_{--}(x,y), \\
b_{+-}(x,y)=a_{+-}(x,y)-a_{+-}(0,0),\quad
b_{-+}(x,y)=a_{-+}(x,y)-a_{-+}(0,0)
\endgather
$$
vanish at $(0,0)$. Each of these four functions can be written in the
form $x(\dots)+y(\dots)$ where the expressions denoted by dots are
certain analytic (hence bounded) functions. Since the expressions
$\frac x{x+y}$ and $\frac y{x+y}$ are bounded, we have to examine
the functions of the form $x^\varepsilon y^\delta$ where
$\varepsilon,\delta$ take the values $\pm\mu$. Since $\mu$ is pure
imaginary or $|\mu|<1/2$, the latter functions are square
integrable. This completes the proof of Claim 2.

Thus, we have proved the proposition for the case $\mu\ne0$. In the
logarithmic case $\mu=0$ the argument is quite similar; we shall only
indicate necessary modifications.

We have from \tht{2.5}
$$
x^{-1/2}W_{\ka,0}=\lim_{\mu\to0}x^{-1/2}W_{\ka,\mu}(x)
=\ln x\cdot a_0(x)+a_1(x),  \tag2.24
$$
where $a_0$ and $a_1$ are certain analytic functions such that
$$
\gathered
a_0(0)=-\frac1{\Ga(\frac12-\ka)}\\
a_1(0)=-\frac1{\Ga(\frac12-\ka)}(\psi(\frac12-\ka)-2\psi(1)),\,
\quad \psi(\cdot):=\frac{\Ga'(\cdot)}{\Ga(\cdot)}\,
\endgathered
\tag2.25
$$
see formulas \tht{6.9(2)} and \tht{6.7(13)} in \cite{E1}.
{}From this we get
$$
C(x,y)=\frac{\ln x\ln y\cdot a_{00}(x,y)+\ln x \cdot a_{01}(x,y)
+\ln y\cdot a_{10}(x,y)+a_{11}(x,y)}{x+y}\,,
$$
where $a_{00}$, $a_{01}$, $a_{10}$, $a_{11}$ are certain analytic
functions. Further, using \tht{2.25} and the well--known relation
$\psi(1+a)-\psi(a)=a^{-1}$ we get
$$
a_{00}(0,0)=0,\quad a_{01}(0,0)=-a_{10}(0,0).
$$
Now we set
$$
C_0(x,y)=\frac{\ln x\cdot a_{01}(0,0)+\ln y\cdot a_{10}(0,0)}{x+y}
=a_{01}(0,0)\,\frac{\ln\left(\frac xy\right)}{x+y}\,,
$$
we define $C'(x,y)$ as above and we state the same two claims as
above, which imply the proposition.

To check Claim 1 we again pass to new variables and transform
the kernel $C(x,y)$ to
$$
const\,\frac{\xi-\eta}{e^{(\xi-\eta)/2}+e^{-(\xi-\eta)/2}}\,.  \tag2.26
$$
Employing the formula
$$
Fourier_u\left\{\frac{v}{e^{v/2}+e^{-v/2}}\right\}
=\pi^2\,\frac{\sh(\pi u)}{\ch^2(\pi u)}  \tag2.27
$$
we get that after the Fourier transform the kernel \tht{2.26} becomes
multiplication by the function \tht{2.27} times a scalar factor. Then we
remark that the latter function is bounded.

As for Claim 2, it is verified in exactly the same way as in the
nonlogarithmic case: here we employ the fact that the function \tht{2.24}
has exponential decay at infinity and is square integrable near zero.

This completes the proof. \qed
\enddemo

\proclaim{Theorem 2.4} Consider the operator $K$ in the Hilbert space
$L^2(\R_+,dx)\oplus L^2(\R_+,dx)$ defined by the matrix Whittaker
kernel \tht{2.11}--\tht{2.12}, and assume that the parameter
$a=\frac{z+z'}2$ satisfies the condition $-1/2<a<1/2$. Set
$$
L=\bmatrix 0 & A \\ -B & 0 \endbmatrix\,, \tag2.28
$$
where
$$
\gather
A(x,y)=D(y,x)=\frac\si\pi\, \left(\frac xy\right)^{-a}
\frac{e^{-\frac{x+y}2}}{x+y}\tag2.29 \\
B(x,y)=D(x,y)=A(y,x)=\frac\si\pi\, \left(\frac xy\right)^a
\frac{e^{-\frac{x+y}2}}{x+y}\,.  \tag2.30
\endgather
$$
{\rm(}Note that $L$ is bounded because $D$ is bounded by Proposition
2.2 {\rm)}.

Then we have
$$
K=\frac L{1+L}\,.  \tag2.31
$$
\endproclaim

\demo{Proof} We know from Propositions 2.2, 2.3 that the kernels
$C(x,y)$, $D(x,y)$ (see \tht{2.14}, \tht{2.15}) define bounded
operators $C$, $D$. Together with Proposition 2.1 this means that in
the formula \tht{2.13} for the matrix Whittaker kernel we may
interpret all the products as products of bounded operators. In
particular, this implies that $K$ is bounded.  

Furthermore, the fact that the matrix Whittaker kernel is
$J$-Hermitian implies the following operator relation:
$$
K=\bmatrix CD & C\\ DCD-D & DC \endbmatrix
=\bmatrix D^*C^* & D^*-D^*C^*D^*\\ -C^* & C^*D^* \endbmatrix\,.
\tag2.32 
$$
(Here the adjoint operators $(\dots)^*$ coincide with the transposed
ones because all the operators are real.)

The desired relation \tht{2.13} is equivalent to $K+KL=L$, which in turn
means that
$$
\bmatrix CD & C\\ DCD-D & DC \endbmatrix +
\bmatrix CD & C\\ DCD-D & DC \endbmatrix \,
\bmatrix 0 & D^*\\ -D& 0\endbmatrix =
\bmatrix 0 & D^*\\ -D& 0\endbmatrix \,,
$$
or
$$
\bmatrix 0 & C+CDD^*\\ -D & DC+DCDD^*-DD^*\endbmatrix =
\bmatrix 0 & D^*\\ -D& 0\endbmatrix\,,
$$
or
$$
C+CDD^*=D^*, \quad DC+DCDD^*-DD^*=0.
$$
Now, the latter two relations are direct consequences of \tht{2.32}.
\qed 
\enddemo

The result seems to be quite surprising. First, the expression for the
kernel $L$ is very simple and involves no special functions.  Second,
let us pass from the parameters $z,z'$ to the parameters $a$,
$\mu$; then we see that in \tht{2.28}--\tht{2.30}, $\mu$ occurs only
in the scalar factor $\si$. 

\proclaim{Corollary 2.5} Let $a$ be fixed, $-1/2<a<1/2$. Then the
operators $K$ corresponding to various values of the parameter $\mu$
pairwise commute. The same holds for the blocks $K_{++}$ or the
blocks $K_{--}$.
\endproclaim

\demo{Proof} According to \tht{2.29}, \tht{2.30} we can write
$$
A=\si A_0,\quad B=\si B_0
$$
where the operators $A_0$ and $B_0$ do not depend on $\mu$.
Consequently, 
$$
L=\si\,\bmatrix 0 & A_0 \\ -B_0 & 0 \endbmatrix
$$
This means that when $\mu$ varies, the operators $L$ differ by a
scalar fact only. Since $K=L(1+L)^{-1}$, we conclude that the
corresponding operators $K$ form a commutative family. 

Next, by \tht{1.15}, 
$$
K_{++}=AB(1+AB)^{-1}.  \tag2.33
$$ 
Since $AB=\si^2 A_0B_0$ where $A_0B_0$ does not depend on $\mu$, the
operators $AB$ with various $\mu$ form a commutative family. So, the
same holds for the operators $K_{++}$.

For the blocks $K_{--}$ the argument is the same. \qed
\enddemo 

This suggests the idea to take as the parameters the couple $a,\si$.

\head 3. Diagonalization of the kernels \endhead

Now we are in a position to perform the
spectral analysis of the kernel $L$ --- to ``diagonalize'' it in a
continual basis and hence to ``diagonalize'' the kernel $K$, too. 

\proclaim{Proposition 3.1} The ``ordinary'' Whittaker kernel $K_{++}$
with parameters $a,\mu$ 
commutes with the Sturm--Liouville differential operator
$$
\D(a)=-\frac{d}{dx}\,x^2\,\frac{d}{dx}\, + \, 
\left(a-\frac x2\right)^2,\tag3.1
$$
i.e., the kernel satisfies the following differential equation
$$
\D(a)_x K_{++}(x,y)=\D(a)_y K_{++}(x,y), \tag3.2
$$
where the subscript $x$ or $y$ indicates the variable on which the
differential operator acts.
\endproclaim 

\demo{Proof} This is a limit case of Proposition 6.2 in \cite{P.III} and
can be verified by a direct computation. We do not give a detailed
proof, because we shall employ this result as a prompt only. \qed
\enddemo 

It is worth noting that $\D(a)$ does not depend on $\mu$, which
agrees with the fact that the kernels $K_{++}$ with
varying $\mu$ form a commutative family (Corollary 2.5). 

Consider the following functions on $\R_+$:
$$
f_{a,m}(x)=\frac 1x \, W_{a,im}(x),\quad m>0. \tag3.3
$$
We have
$$
\D(a)f_{a,m}=(a^2+\frac14+m^2)f_{a,m}. \tag3.4
$$
According to \cite{W}, the functions
$f_{a,m}$ with $a$ fixed and $m$ ranging over $\R_+$ form a continual
basis in $L^2(\R_+)$ diagonalizing $\D(a)$. Moreover, an explicit
Plancherel formula holds:
$$
(f,g)_{L^2(\R_+)}=\int_0^{+\infty}
\frac{(f,f_{a,m})(f_{a,m},g)}{(f_{a,m},f_{a.m})}\,dm\,, \tag3.5
$$
where 
$$
(f_{a,m},f_{a.m})
:=\frac{\pi^2}{\Ga(\frac12-a-im)\Ga(\frac12-a+im)}\,. \tag3.6
$$

Consider the decomposition \tht{2.1} and take in its first component the
basis $\{f_{a,m}\}_{m\ge0}$ and in its second component --- the basis
$\{f_{-a,m}\}_{m\ge0}$. Together they form a (continual) basis in the
whole space $L^2(\R^*)$. The following claim describes the
diagonalization of $L$ in this basis.

\proclaim{Theorem 3.2} Let $A$ and $B$ be as in \tht{2.29}, \tht{2.30}, and assume that $|a|<\frac12$. Then we have
$$
\gather
Af_{-a,m}=\frac\si\pi\, 
\Ga(\frac12-a+im)\Ga(\frac12-a-im)f_{a,m}\tag3.7a\\
Bf_{a,m}=\frac\si\pi\, 
\Ga(\frac12+a+im)\Ga(\frac12+a-im)f_{-a,m}. \tag3.7b
\endgather
$$
\endproclaim

\demo{Proof} The function $Af_{-a,m}$ is essentially the Stieltjes
transform of the function 
$$
y^{-a-1}\exp(-\frac12 y)W_{-a,im}(y),
$$
which is given, under the assumption $\Re(-a)>-\frac12$, in
\cite{E2, 14.3(53)}. The expression for $B=A^*$ is
obtained in the same way; here we need 
$\Re(a)>-\frac12$. Thus, we have entirely used the assumption
$|a|<\frac12$.   \qed
\enddemo

\example{Remark 3.3} We have
$$
ABf_{a,m}=\frac{\cos(2\pi\mu)-\cos(2\pi a)}
{\cos(2\pi im)+\cos(2\pi a)}\, f_{a,m}\,. \tag3.8
$$
This shows that $a=\pm\frac12$ are ``critical'' points: when
$|a|<\frac12$, the operator $AB=AA^*$ is bounded as its spectrum is
bounded, whence $A$ and $B$ are bounded. But in the limit
$a\to\pm\frac12$ the spectrum of $AB$ becomes unbounded. 
\endexample

Employing \tht{2.33}, we get, as a corollary of \tht{2.18}, a
diagonalization of $K_{++}$: 

\proclaim{Corollary 3.4} Assume $|a|<\frac12$. We have
$$
K_{++}f_{a,m}=\frac{\cos(2\pi\mu)-\cos(2\pi a)}
{\cos(2\pi\mu)+\cos(2\pi i m)}\, f_{a,m}\,.
$$
\endproclaim
\example{Remark 3.5}
The operator $A$ for $a=0$ arose before in the asymptotic analysis of the
Painlev\'e transcendent of the third kind [MTW] and two--dimensional Ising model [T]. It has also been used in approximating the resolvent of a certain more difficult operator in [TW4].

The spectral analysis of the kernel (2.29) for $a=0$ was carried out in [MTW], see also [T].

If $a=0$ then (3.3), (3.7a) turns into 
$$
f_{0,m}=\frac 1x W_{0,im}(x)=\frac{1}{\sqrt{\pi x}}\, K_{im}\left(\frac x2\right),
$$ 
$$
Af_{0,m}=\frac{\sigma}{\ch \pi m}\,f_{0,m}
$$
which agrees with the results of [MTW] (here $K_{\nu}(x)$ stands for the Bessel $K$--function).
\endexample
\head 4. The matrix tail kernel \endhead

In \cite{P.III}, we studied a `tail kernel' associated to the Whittaker
kernel. This is a translation invariant kernel on $\Bbb R$
generalizing the sine kernel. Specifically, it has the form
$$
\frac BA\,\frac{\sin(A(\xi-\eta))}{\sh(B(\xi-\eta))}
\quad \text{or}\quad
\frac{B\sh(A(\xi-\eta))}{A\sh(B(\xi-\eta))} \tag4.1
$$
depending on whether $\mu$ is pure imaginary or real, respectively; $A$ and
$B$ are certain constants depending on $z,z'$ (see also below). Here
we aim to study a similar object for the {\it matrix\/} Whittaker
kernel. 

The proofs are omitted; they are quite similar to that given in
\cite{P.III}. For the sake of simplicity, we shall assume $z\ne z'$. 

Let us briefly recall how the tail kernel arises. According to
\cite{P.II}, Theorem 4.1.1, the Whittaker kernel $K_{++}(x,y)$
behaves near 
$(0,0)$ as follows
$$
K_{++}(x,y)\approx\frac{\sin(\pi z)\sin(\pi z')}
{\pi\sin(\pi(z-z'))}\,
\frac1{\sqrt{xy}}\,
\frac{(x/y)^{\frac{z-z'}2}-(x/y)^{-\frac{z-z'}2}}
{(x/y)^{\frac12}-(x/y)^{-\frac12}}  \tag4.2
$$
In particular, on the diagonal,
$$
K_{++}(x,x)\approx\frac Cx\,,
$$
where
$$
C=C(z,z')=\frac{(z-z')\sin(\pi z)\sin(\pi z')}
{\pi\sin(\pi(z-z'))}\,.  \tag4.3
$$

We pass to new variables to make the density function $K_{++}(x,x)$
asymptotically equal to 1. Specifically, we take
$$
x=e^{-\xi/C}, \quad y=e^{-\eta/C}\,.
$$
Then the resulting kernel in $\xi,\eta$ takes the form
$$
\Kpp(\xi,\eta)\,+\,\text{remainder term}
$$
where the remainder term tends to zero as $\xi,\eta\to+\infty$ and $\Kpp$
is a translation invariant kernel equal to 1 on the diagonal,
$$
\Kpp(\xi,\eta)=\frac1{z-z'}\,\frac{\sh(A(\xi-\eta))}
{\sh(B(\xi-\eta))}  \tag4.4
$$
with
$$
\gather
B=\frac1{2C}=\frac{\pi\sin(\pi(z-z'))}
{2(z-z')\sin(\pi z)\sin(\pi z')}  \tag4.5\\
A=(z-z')B=\frac{\pi\sin(\pi(z-z'))}
{2\sin(\pi z)\sin(\pi z')}\,.  \tag4.6
\endgather
$$
Note that $B$ is real and strictly positive while $A$ is real or pure
imaginary together with $\mu$ 
(so, for pure imaginary $\mu$ the kernel actually has the form of the
first expression in \tht{4.1}).

\proclaim{Theorem 4.1} Application of the same procedure to the
matrix Whittaker kernel 
$$
K(x,y)=\bmatrix K_{++}(x,y) & K_{+-}(x,y)\\
K_{-+}(x,y) & K_{--}(x,y)\endbmatrix  \tag4.7
$$
leads to a translation invariant block kernel in $\xi, \eta$,
$$
\Cal K(\xi,\eta)=\bmatrix \Kpp(\xi,\eta) & \Kpm(\xi,\eta) \\
\Kmp(\xi,\eta) & \Kmm(\xi,\eta) \endbmatrix\,,  \tag4.8
$$
where all the blocks are real, 
$$
\Kpp(\xi,\eta)=\Kmm(\xi,\eta), \quad \Kpm(\xi,\eta)=-\Kmp(\eta,\xi),
\tag4.9 
$$
$\Kpp$ is as in \tht{4.4} and
$$
\Kpm(\xi,\eta)=\frac1{\sqrt{\sin(\pi z)\sin(\pi z')}}\,
\frac1{z-z'}\,
\frac{\sin(\pi z)e^{A(\xi-\eta)}-\sin(\pi z')e^{-A(\xi-\eta)}}
{e^{B(\xi-\eta)}+e^{-B(\xi-\eta)}}  \tag4.10
$$
with $A$ and $B$ as in \tht{4.5}, \tht{4.6}. 

I.e., in the new variables, the kernel \tht{4.7} takes the form
\tht{4.8} plus a remainder term which tends to zero as $\xi,\eta\to+\infty$.
\endproclaim

\demo{Idea of proof} This claim is a generalization of Theorem 3.2
from  \cite{P.III} and is proved in the same way. In addition
to \tht{4.2} we employ an asymptotic formula for the kernel
$K_{+-}(x,y)$ near $(0,0)$:
$$
K_{+-}(x,y)\approx\frac{\sqrt{\sin(\pi z)\sin(\pi z')}}
{\pi\sin(\pi(z-z'))}\,
\frac1{\sqrt{xy}}\,
\frac{\sin(\pi z)(x/y)^{\frac{z-z'}2}+\sin(\pi z')(x/y)^{-\frac{z-z'}2}}
{(x/y)^{\frac12}+(x/y)^{-\frac12}}\,,  \tag4.11
$$
which is proved similarly. \qed
\enddemo

Consider the integral operator in the Hilbert space of square integrable
$\C^2$-valued functions on $\R$ that is defined by the kernel 
$\Cal K(\xi, \eta)$. Since the kernel is translation invariant, the
operator in question is a convolution operator. Under the Fourier
transform it turns into the operator of multiplication by a 
$2\times2$ matrix--valued function, say, $\KK(u)$.

\proclaim{Proposition 4.2} The above defined matrix function on $\R$
has the form 
$$
\KK(u)=\bmatrix f(u) & g(u)\\ -\bar g(u) & f(u) \endbmatrix\,,  \tag4.12
$$
where $f(u)$ is a real function, $g(u)$ is a complex function, 
$\bar g(u)=\overline{g(u)}$, 
$$
\gather
f(u)=2\sin(\pi z)\sin(\pi z')\,
\frac1{\cos(\pi(z-z'))+\ch(\pi u/B)}\\
g(u)=2\sqrt{\sin(\pi z)\sin(\pi z')}\,
\frac{\cos(\pi(z+z')/2+i\pi u/(2B))}
{\cos(\pi(z-z'))+\ch(\pi u/B)}\,.
\endgather
$$
\endproclaim

\demo{Sketch of proof} Rewrite the expressions \tht{4.4}, \tht{4.10}
in the form
$$
\Kpp(\xi,\eta)=\kpp(\xi-\eta), \quad
\Kpm(\xi,\eta)=\kpm(\xi-\eta)\,
$$
where $\kpp$ and $\kpm$ are functions of a single variable, say,
$\zeta$. It follows from the symmetry properties \tht{4.9}
that $\KK(u)$ has the form \tht{4.12}, where $f$ and $g$ are the Fourier
transforms of $\kpp$ and $\kpm$, respectively: 
$$
f(u)=\int e^{i\pi u\zeta}\,\kpp(\zeta)d\zeta, \quad
g(u)=\int e^{i\pi u\zeta}\,\kpm(\zeta)d\zeta\,.
$$
Since $\kpp$ is an even function, $f(u)$ is real.

The desired explicit expression for $f(u)$ is a table integral, see,
e.g., \cite{E2, 1.9(14)}. The explicit expression for
$g(u)$ can be derived from another table integral, see 
\cite{E2, 3.2(15)}. \qed 
\enddemo

\head 5. Degeneration to a Bessel--type
kernel \endhead

We shall need three Bessel functions: the Bessel function of the
first kind
$$
J_\nu(X)=\frac{(X/2)^\nu}{\Ga(\nu+1)}\,
\F(\nu+1;-(X/2)^2)=\frac{(X/2)^\nu}{\Ga(\nu+1)}
\sum_{m\ge0}\frac{(-1)^m(X/2)^{2m}}{m!(\nu+1)_m}\,,
$$
the modified Bessel function of the first kind
$$
I_\nu(X)=\frac{(X/2)^\nu}{\Ga(\nu+1)}\,
\F(\nu+1;(X/2)^2)=\frac{(X/2)^\nu}{\Ga(\nu+1)}
\sum_{m\ge0}\frac{(X/2)^{2m}}{m!(\nu+1)_m}\,,
$$
and the modified Bessel function of the third kind, also called the
Macdonald function
$$
K_\nu(X)=\frac{\pi}{2\sin(\pi\nu)}(I_{-\nu}(X)-I_\nu(X)).
$$
Here we assume $X>0$.

We fix the parameters $z_0$, $z'_0$ satisfying the assumptions
\tht{2.2}, and we set 
$$
a_0=\frac{z_0+z'_0}2\,,\quad \mu=\frac{z_0-z'_0}2\,.
$$ 
In terms of $a_0, \mu$, the restrictions on $z_0,z'_0$ take the
form \tht{2.4}.

Let $N$ be an integer; then the parameters 
$$
z=z_0+N, \quad z'=z'_0+N
$$
will also satisfy the same restrictions \tht{2.2} as $z_0$, $z'_0$. We set
$$
a:=\frac{z+z'}2
$$
and note that 
$$
a=a_0+N\,, \quad \frac{z-z'}2=\mu.
$$

We introduce four functions in a positive variable $\xi$:
$$
\gather
A(\xi)=\frac{\sin(\pi z_0)J_{2\mu}(2\sqrt{\xi})-
\sin(\pi z'_0)J_{-2\mu}(2\sqrt{\xi})}{\sin(2\pi\mu)}\,,\\
B(\xi)=K_{2\mu}(2\sqrt{\xi})\,,\\
\wtA(\xi)=\sqrt{\xi}\,\frac{\sin(\pi z_0)J'_{2\mu}(2\sqrt{\xi})-
\sin(\pi z'_0)J'_{-2\mu}(2\sqrt{\xi})}{\sin(2\pi\mu)}\,,\\
\wtB(\xi)=\sqrt{\xi}\,K'_{2\mu}(2\sqrt{\xi})\,,
\endgather
$$
where 
$$
J'_\nu(X)=\frac{d}{dX}J_\nu(X),\quad
K'_\nu(X)=\frac{d}{dX}K_\nu(X).
$$

Note that
$$
\wtA(\xi)=\left(\xi\frac{d}{d\xi}\right)A(\xi), \quad
\wtB(\xi)=\left(\xi\frac{d}{d\xi}\right)B(\xi),
$$
because, for a function $f(\cdot)$, 
$$
\left(\xi\frac{d}{d\xi}\right)f(2\sqrt{\xi})=\sqrt{\xi}f'(2\sqrt{\xi}).
$$

Finally, we let $N\to\infty$ and associate with the positive
variables $x,y$ the `scaled variables' $\xi,\eta$,
$$
x=\xi/N,\quad y=\eta/N.
$$

\proclaim{Theorem 5.1} In the scaled limit, as $N\to+\infty$ inside
$2\Bbb Z$, the matrix Whittaker kernel in the variables $x,y$,
$$
K(x,y)=\bmatrix K_{++}(x,y) & K_{+-}(x,y)\\
K_{-+}(x,y) & K_{--}(x,y)\endbmatrix
$$
tends to a matrix kernel $K^{lim}$ in the variables $\xi,\eta$ with the
following blocks:
$$
\gather
K^{lim}_{++}(\xi,\eta)
=\frac{A(\xi)\wtA(\eta)-\wtA(\xi)A(\eta)}{\xi-\eta}\,,\tag5.1a \\
K^{lim}_{+-}(\xi,\eta)=-\,\frac{2\sqrt{\sin(\pi z_0)\sin(\pi z'_0)}}{\pi}\,
\frac{A(\xi)\wtB(\eta)-\wtA(\xi)B(\eta)}{\xi+\eta}\,,\tag5.1b \\
K^{lim}_{-+}(\xi,\eta)=-K^{lim}_{+-}(\eta,\xi),\tag5.1c \\
K^{lim}_{--}(\xi,\eta)=\frac{4\sin(\pi z_0)\sin(\pi z'_0)}{\pi^2}\,
\frac{B(\xi)\wtB(\eta)-\wtB(\xi)B(\eta)}{\xi-\eta}\,.\tag5.1d
\endgather
$$
\endproclaim

\example{Comments} 1) Under the shift $z_0\mapsto z_0+1$,
$z'_0\mapsto z'_0+1$, the functions $B$, $\wtB$ remain stable while
the functions $A$, $\wtA$ are multiplied by $-1$. It follows that
under this shift, the diagonal blocks $K^{lim}_{++}$, $K^{lim}_{--}$
are stable while the blocks $K^{lim}_{+-}$, $K^{lim}_{-+}$ are
multiplied by $-1$. We could assume $N$ tends to infinity inside
$\Bbb Z$ (instead of $2\Bbb Z$) by introducing in the blocks
$K_{+-}$, $K_{-+}$ the extra factor $(-1)^{N}$; such a factor does
not affect the correlation functions.

2) Note that the function $B$ depends only on
the parameter $\mu$ while $A$ depends on the both parameters
$a_0,\mu$. This results in a strong asymmetry between the diagonal
blocks $K^{lim}_{++}$ and $K^{lim}_{--}$. Of course, these blocks are
interchanged if we let $N$ tend to $-\infty$ instead of $+\infty$. 

3) When $\mu$ is real and one of the parameters $z_0=a_0+\mu$,
$z'_0=a_0-\mu$ becomes integer, the kernel $K_{++}$ degenerates to
the conventional Bessel kernel \footnote{About the Bessel kernel, see
\cite{F, NS, NW, TW2}.} 
$$
\frac{J_\nu(2\sqrt{\xi})\sqrt{\eta}J'_\nu(2\sqrt{\eta})
-\sqrt{\xi}J'_\nu(2\sqrt{\xi})J_\nu(2\sqrt{\eta})}
{\xi-\eta}\,, \quad \nu=2\mu. \tag5.2
$$
This agrees with the degeneration of the Whittaker kernel to the
Laguerre kernel, see \cite{P.III}, Remark 2.4. \footnote{It is well known
that the Bessel kernel can be obtained in a scaling limit of the
Laguerre kernel, see \cite{F, NS, NW, TW2}.} Thus, thanks to the
parameter $a_0$, the kernel $K^{lim}_{++}$ provides a deformation of
the Bessel kernel; one more new point is that the index $\mu$ in the
expression for $K^{lim}_{++}$ can be pure imaginary.

4) The kernel 
$$
K^{lim}_{--}(\xi,\eta)=const\, 
\frac{K_{2\mu}(2\sqrt{\xi})\sqrt{\eta}K'_{2\mu}(2\sqrt{\eta})
-\sqrt{\xi}K'_{2\mu}(2\sqrt{\xi})K_{2\mu}(2\sqrt{\eta})}
{\xi-\eta} \tag5.3
$$
except the scalar factor $const$, depends only on
$\mu$ and looks quite similar to the Bessel kernel \tht{5.2}. One could
call it the {\it Macdonald kernel}. 

5) The scalar factor $const$ in \tht{5.3} can
be written in the form
$$
const=\frac2{\pi^2}\,(\cos(2\pi\mu)-\cos(2\pi a_0)). \tag5.4
$$
This expression is periodic in $a_0$ with period 1. It is strictly
positive (because of the assumptions on the parameters). 
When $\mu$ is fixed, its maximal value, attained at the point
$a_0=1/2$, is equal to $4\cos^2(\pi \mu)/\pi^2$. 
\endexample

\demo{Proof} {\it Step 1: A transformation of the matrix Whittaker
kernel.\/} Recall the expression of the Whittaker function through
the confluent hypergeometric function:
$$
\aligned
x^{-1/2}e^{x/2}W_{\ka,\mu}(x) &=
\frac{\Ga(-2\mu)x^\mu}{\Ga(\frac12-\ka-\mu)}\,
\FF(\tfrac12-\ka+\mu;2\mu+1;x)\\
&+\frac{\Ga(2\mu)x^{-\mu}}{\Ga(\frac12-\ka+\mu)}\,
\FF(\tfrac12-\ka-\mu;-2\mu+1;x).
\endaligned \tag5.5
$$

Let us abbreviate
$$
\gather
\varphi(x)=x^{-1/2}e^{x/2}W_{a+\frac12,\mu}(x),\tag5.6a \\
\varphi_-(x)=x^{-1/2}e^{x/2}W_{a-\frac12,\mu}(x),\tag5.6b \\
\psi(x)=x^{-1/2}e^{x/2}W_{-a+\frac12,\mu}(x),\tag5.6c \\
\psi_-(x)=x^{-1/2}e^{x/2}W_{-a-\frac12,\mu}(x),\tag5.6d \\
\wtphi(x)=x\varphi'(x),\\
\wtpsi(x)=x\psi'(x).
\endgather
$$

{}From \tht{5.5}, \tht{5.6a,b,c,d} and the series expansion
$$
\FF(\al;\ga;x)=\sum_{m\ge0}\frac{(\al)_m}{m!(\ga)_m}x^m
$$
we readily get
$$
\gather
\varphi_-(x)=\frac{\wtphi(x)-a\varphi(x)}{zz'}\,,\\
\psi_-(x)=\frac{\wtpsi(x)+a\psi(x)}{zz'}\,.
\endgather
$$

Together with the definition of the matrix Whittaker kernel this
implies 
$$
\align
e^{\frac{x+y}2}K_{++}(x,y)&=
\frac1{\Ga(z)\Ga(z')}\,
\frac{\varphi(x)\varphi_-(y)-\varphi_-(x)\varphi(y)}{x-y}\\
&=\frac1{\Ga(z)\Ga(z')zz'}\,
\frac{\varphi(x)\wtphi(y)-\wtphi(x)\varphi(y)}{x-y}\,, \tag5.7a
\endalign
$$
$$
\align
e^{\frac{x+y}2}K_{--}(x,y)&=
\frac1{\Ga(-z)\Ga(-z')}\,
\frac{\psi(x)\psi_-(y)-\psi_-(x)\psi(y)}{x-y}\\
&=\frac1{\Ga(-z)\Ga(-z')zz'}\,
\frac{\psi(x)\wtpsi(y)-\wtpsi(x)\psi(y)}{x-y}\,, \tag5.7b
\endalign
$$
$$
\gather
e^{\frac{x+y}2}K_{+-}(x,y)=
\frac{\sqrt{\sin(\pi z)\sin(\pi z')}}{\pi}\,
\frac{\varphi(x)\psi(y)+zz'\varphi_-(x)\psi_-(y)}{x+y}\\
=\frac{\sqrt{\sin(\pi z)\sin(\pi z')}}{\pi}\\
\times\,\frac{-\frac{a}{zz'}(\varphi(x)\wtpsi(y)-\wtphi(x)\psi(y))
+\left(1-\frac{a^2}{zz'}\right)\varphi(x)\psi(y)+
\frac1{zz'}\wtphi(x)\wtpsi(y)}{x+y}\,. \tag5.7c
\endgather
$$

{\it Step 2: The scaling limit of the functions $\varphi$, $\wtphi$,
$\psi$, $\wtpsi$.\/} 
We start with the well--known limit formula, which is readily
obtained from the standard series expansions for $\FF$ and $\F$:
$$
\lim_{|\al|\to\infty}\FF\left(\al;\ga;\frac{\xi}{\al}\right)=
\F(\ga;\xi). \tag5.8
$$
Here $\al$, $\ga$, $\xi$ are allowed to be any complex numbers with
the only restriction $\ga\ne0,-1,-2,\dots$. The convergence is
uniform on compact sets in the $\xi$-plane, which implies that this
limit relation can be differentiated with respect to $\xi$.

It follows that
$$
\aligned
&\phantom{\sim}\frac{\Ga(-2\mu)x^\mu}{\Ga(-a-\mu)}\,
\FF(-a+\mu;2\mu+1;x)\\
&\sim\,\frac{\Ga(-2\mu)\xi^\mu}{\Ga(-a-\mu)N^\mu}\,\F(2\mu+1;-\xi)\\
&=\,\frac{\Ga(-2\mu)\Ga(2\mu+1)}{\Ga(-a-\mu)N^\mu}\,
\frac{\xi^\mu}{\Ga(2\mu+1)}\,\F(2\mu+1;-\xi)\\
&\sim\,\Ga(a+1)\,\frac{\sin(\pi z)}{\sin(2\pi\mu)}\,
J_{2\mu}(2\sqrt{\xi}).
\endaligned \tag5.9
$$
Here we have used the relations $z=a+\mu$, $a\sim N$, 
$$
\Ga(-w)\Ga(1+w)=-\,\frac{\pi}{\sin(\pi w)}\,,\quad
\frac{\Ga(N+const_1)}{\Ga(N+const_2)}\,\sim\,N^{const_1-const_2}
$$
and the expression of the Bessel function $J_\nu$ through the $\F$
function. 

Likewise, 
$$
\aligned
&\phantom{\sim}\frac{\Ga(-2\mu)x^\mu}{\Ga(a-\mu)}\,
\FF(a+\mu;2\mu+1;x)\\
&\sim\,\frac{\Ga(-2\mu)\Ga(2\mu+1)}{\Ga(a-\mu)N^\mu}\,
\frac{\xi^\mu}{\Ga(2\mu+1)}\,\F(2\mu+1;\xi)\\
&\sim\,-\,\frac1{\Ga(a)}\,\frac{\pi}{\sin(2\pi\mu)}\,
I_{2\mu}(2\sqrt{\xi}).
\endaligned \tag5.10
$$

By the definition \tht{5.6a} of the function $\varphi(x)$ and the expression
\tht{5.5} for the Whittaker function, $\varphi(x)$ is equal to the
left--hand side of \tht{5.9} plus the symmetric expression obtained by
inserting $-\mu$ in place of $\mu$. Then it follows from \tht{5.9} that
$$
\aligned
\varphi(x)\,&\sim\,\Ga(a+1)\left(
\frac{\sin(\pi z)}{\sin(2\pi\mu)}\,J_{2\mu}(2\sqrt{\xi})+
\frac{\sin(\pi z')}{\sin(-2\pi\mu)}\,J_{-2\mu}(2\sqrt{\xi})
\right)\\
&=\,\Ga(a+1)A(\xi).
\endaligned \tag5.11a
$$
Here we have used the definition of $A(\xi)$ and the fact that
$$
a-\mu=z',\quad \sin(\pi z)=\sin(\pi z_0+\pi N)=\sin(\pi z_0),\quad
\sin(\pi z')=\sin(\pi z'_0)
$$
because $N\in 2\Bbb Z$ by assumption.

Likewise, it follows from \tht{5.6b}, \tht{5.5} and \tht{5.10} that 
$$
\aligned
\psi(x)\,&\sim\,-\,\frac1{\Ga(a)}\left(
\frac{\pi}{\sin(2\pi\mu)}\,I_{2\mu}(2\sqrt{\xi})+
\frac{\pi}{\sin(-2\pi\mu)}\,J_{-2\mu}(2\sqrt{\xi})
\right)\\
&=\frac2{\Ga(a)}\,K_{2\mu}(2\sqrt{\xi})=
\frac2{\Ga(a)}B(\xi).
\endaligned \tag5.11b
$$

Recall that our asymptotic formulas, which are based on the limit formula
\tht{5.8}, are stable under differentiation and note that the
differential operator $x\frac{d}{dx}$ is invariant 
relative to the change of a variable $x\mapsto \xi=xN$. It follows that
$$
\gather
\wtphi(x)=\left(x\frac{d}{dx}\right)\varphi(x)\,
\sim\,\Ga(a+1)\left(\xi\frac{d}{d\xi}\right)A(\xi)
=\Ga(a+1)\wtA(\xi), \tag5.11c\\
\wtpsi(x)=\left(x\frac{d}{dx}\right)\psi(x)\,
\sim\,\frac{2}{\Ga(a)}\left(\xi\frac{d}{d\xi}\right)B(\xi)
=\frac2{\Ga(a)}\wtB(\xi). \tag5.11d
\endgather
$$

{\it Step 3: The scaling limit of the matrix Whittaker kernel.\/} 
It remains to combine the formulas \tht{5.7a,b,c} with the asymptotic
expressions \tht{5.11a,b,c,d}. 

First of all, note that the transformation of a kernel in $x,y$
under a scaling involves the transformation of a differential, say,
$dy$. We have 
$$
\frac{dy}{x\pm y}\,=\,\frac{d\eta}{\xi\pm\eta}\,,
$$
so that in the scaling limit, the denominator $x\pm y$ simply turns
into $\xi\pm\eta$. Next, in the scaling limit, the exponential
factor in the left--hand side of the formulas \tht{5.7a,b,c} is negligible. 

Inserting the asymptotic expressions for $\varphi$ and $\wtphi$ into
\tht{5.7a} and using the relation
$$
\frac{\Ga(a+1)\Ga(a+1)}{\Ga(z)\Ga(z')zz'}\,\sim\,1
$$
we get the desired formula \tht{5.1a}.

Likewise, inserting the asymptotic formulas for $\psi$ and $\wtpsi$
into \tht{5.7b} and using the relation
$$
\aligned
\frac4{\Ga(-z)\Ga(-z')zz'\Ga(a)\Ga(a)} \,
&\sim\,\frac{4\sin(\pi z)\sin(\pi z')}{\pi^2}\\
&=\,\frac{4\sin(\pi z_0)\sin(\pi z'_0)}{\pi^2}
\endaligned
$$
we get the formula \tht{5.1d}.

Now, let us examine the numerator in \tht{5.7c},
$$
-\tfrac{a}{zz'}[\varphi(x)\wtpsi(y)-\wtphi(x)\psi(y)]
+\left(1-\tfrac{a^2}{zz'}\right)[\varphi(x)\psi(y)]+
\tfrac1{zz'}[\wtphi(x)\wtpsi(y)]. \tag5.12
$$
Note that each of the functions $\varphi(\cdot)$, $\wtphi(\cdot)$ is
asymptotically equivalent to a function in a scaled variable times
the factor $\Ga(a+1)$, while each of the function $\psi(\cdot)$,
$\wtpsi(\cdot)$ is asymptotically equivalent to a function in a
scaled variable times the factor $1/\Ga(a)$. It follows that each of
the three expressions in the squared brackets behaves as a function
in the scaled variables $\xi,\eta$ times the factor
$\Ga(a+1)/\Ga(a)\sim N$. 

Further, the coefficients behave as follows
$$
\tfrac a{zz'} =\tfrac a{a^2-\mu^2}\sim N^{-1},\quad
1-\tfrac{a^2}{zz'}=1-\tfrac {a^2}{a^2-\mu^2}=O(N^{-2}),\quad
\tfrac1{zz'}=\tfrac1{a^2-\mu^2}=O(N^{-2}).
$$

This implies that the second and the third summands in \tht{5.12} are
asymptotically negligible. It is readily verified that the
contribution of the first summand yields the desired formula \tht{5.1b}.

Finally, the relation \tht{5.1c} is immediate from the similar
relation between the blocks of the matrix Whittaker kernel.

This concludes the proof. \qed
\enddemo

\enddocument